\documentclass[letterpaper, 10 pt, conference]{ieeeconf}  

\IEEEoverridecommandlockouts                              

\usepackage{color}
\usepackage{graphicx} 
\usepackage{amsmath} 
\usepackage{amssymb}  
\newtheorem{mydef}{Definition}
\newtheorem{mythm}{Theorem}

\newtheorem{rmk}{Remark}

\title{\LARGE \bf
Controllability of a bent 3-link magnetic microswimmer
}

\author{Laetitia Giraldi$^{*}$, Pierre Lissy$^{\dagger}$, Cl\'{e}ment Moreau$^{\ddagger}$, Jean-Baptiste Pomet$^{*}$
\thanks{$^{*}$Team McTAO, INRIA Sophia Antipolis, France}%
\thanks{$^{\dagger}$CEREMADE, Universit\'{e} Paris-Dauphine, Paris, France}%
\thanks{$^{\ddagger}$ENS de Cachan, France}%
}

\begin{document}

\maketitle
\thispagestyle{empty}
\pagestyle{empty}

\begin{abstract}

In this paper, we focus on a variant of a 3-link magnetic microswimmer which consists of three rigid magnetized segments connected by two torsional springs. In particular, we assume that one of the springs is twisted so that the swimmer is not aligned at rest. By acting on it with an external magnetic field, the swimmer twists and moves through the surrounding fluid. By considering the external magnetic field as a control function, we state a local partial controllability result around the equilibrium states. Then, we propose a constructive method to find the magnetic field that allows the swimmer to move along a prescribed trajectory. Finally, we show numerical simulations in which the swimmer moves along a prescribed path.

\end{abstract}

\section{INTRODUCTION}
At a microscopic scale, swimming in water or another similar fluid is a very different matter from the macroscopic one. Indeed, micro-swimmers face a very small Reynolds number (around $10^{-6}$), which means that the intensity of inertial forces is negligible towards those of viscous ones. Due to promising perspectives of medical micro-robots performing delicate tasks inside the human body, interest in the study of micro-swimmers has been recently growing. 

The shapes and propulsion techniques of these new robots could be inspired by biology, since micro-organisms such as sperm cells or bacterias developed efficient ways to move through a surrounding fluid (see \cite{PeyerZhang05}). One direction of research is to use chemical reactions inside the micro-robot to drive it (see \cite{MirkovicZacharia10}). Another technique consists in using an external magnetic field to drive a magnetized swimmer (see \cite{GaoKagan12,Bibette,GhoshFischer09}).  

In this paper, we focus on this type of propulsion, applied on a simple model of micro-swimmer consisting on three magnetized segments linked by elastic joints. Since the  swimmer is supposed to be small, the hydrodynamic interaction between the swimmer
and the fluid can be modeled by the local drag approximation of Resistive Force Theory introduced in \cite{GrayHancock55}. Such models, with different numbers of segments, have been studied for instance in \cite{Or14} and \cite{alouges2013self}, in which the authors show that sinusoidal magnetic fields allow the swimmer to move forward in a prescribed direction.
 
In \cite{GiraldiPomet16}, the authors show a local controllability result for the 2-segment model around its straight position. 
In this paper, we focus on a 3-segment magnetized micro-swimmer, under the assumption that it is not aligned at its equilibrium. By considering the external magnetic field as a control function, we study how to control the position of the swimmer without prescribing any constraints on the orientation and shape of the swimmer, i.e., we state a local \textit{partial} controllability result.
Then, we develop a constructive method to find a magnetic field such that the robot can move along prescribed paths.

The paper is organized as follows. In Section \ref{model}, we detail the dynamics of the model and state the main result of local \emph{partial} controllability for the bent swimmer. In Section \ref{explicit}, we describe a  practical method that explicitly compute the magnetic field to make the swimmer follow some prescribed trajectory, as soon as the swimmer does not go through its aligned position. Using the latter procedure, we give in Section \ref{numericalsimu} some numerical simulations with leads to control the swimmer along some prescribed trajectories. Finally, Section \ref{perspectives} is dedicated to  some perspectives of this work.

\section{MICROSWIMMER MODEL AND CONTROLLABILITY ISSUES \label{model}}

\subsection{Formulation of the Problem}

\begin{figure}
\begin{center}
\includegraphics[width=3.2in]{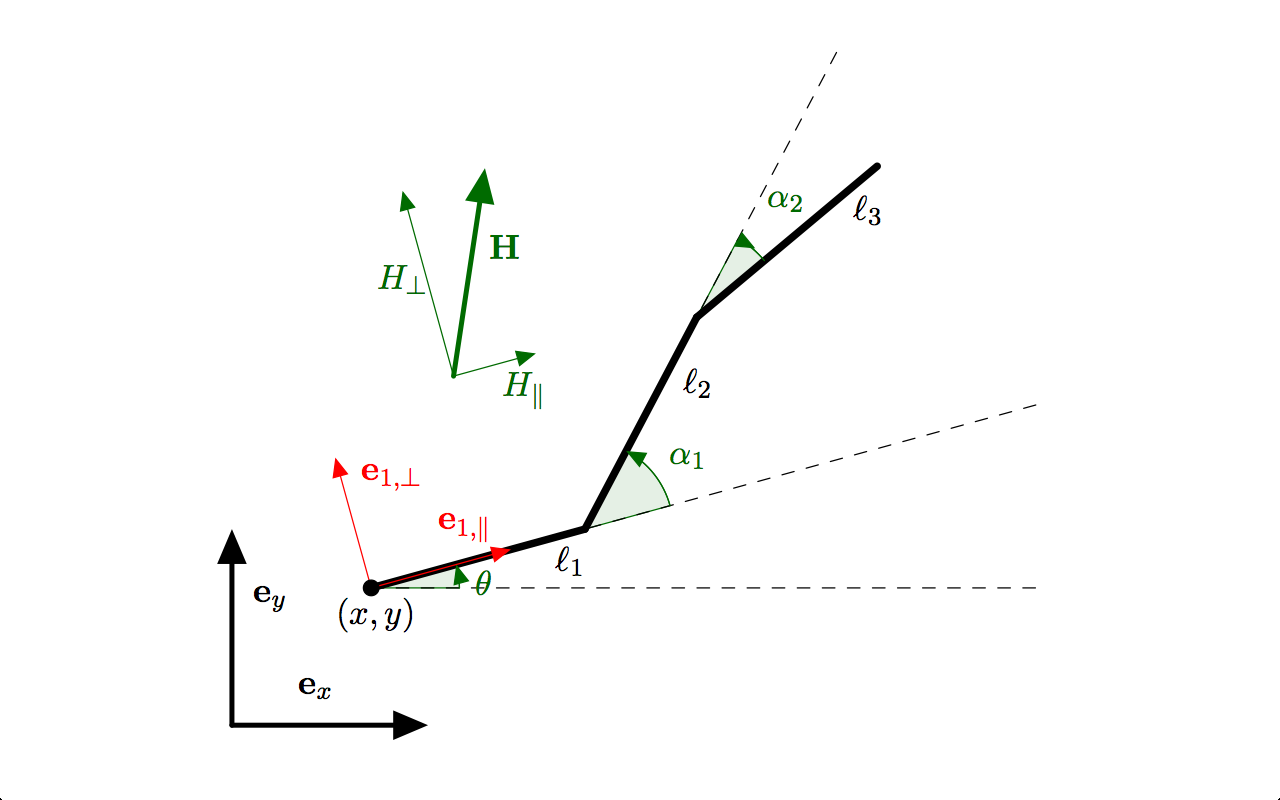}
\caption{Model used for the 3-link microswimmer.}
\label{schema_nageur}
\end{center}
\end{figure}

 We follow the notations, assumption and modelisation introduced in \cite{alouges2015can} and \cite{GiraldiPomet16}. In the present paper, we focus on a micro-swimmer consisting on 3 rigid magnetized segments connected by two torsional springs with stiffness $\kappa$, subject to an external uniform in space magnetic field $\mathbf{H}$. The 3 segments, called $S_1$,$S_2$ and $S_3$, have same length $\ell$, same hydrodynamic drag coefficients $\xi$ and $\eta$, and respective magnetic moments $M_1$, $M_2$ and $M_3$.  The choice of the numerical values for the parameters will be  detailed later in Section \ref{numericalsimu} (see Table \ref{table_values}). The swimmer can move in the 2d-plane defined by the vectors $\mathbf{e}_x$ and $\mathbf{e}_y$. Let us define $\mathbf{e}_z=\mathbf{e}_x \times \mathbf{e}_y$. Let $\mathbf{x}=(x,y)$ be the coordinates of the end of $S_1$, we call $\theta$ the angle between $(Ox)$ and $S_1$, and $\alpha_1$ and $\alpha_2$ the angles between $S_1$ and $S_2$ and between $S_2$ and $S_3$. The swimmer is then completely described by 5 variables $(x,y,\theta,\alpha_1,\alpha_2)$ where the pair $(x,y)$ represents the position of the swimmer, $\theta$ its orientation and the pair $(\alpha_1,\alpha_2)$ its shape.  Let us also define the moving frames associated to $S_i$ for $i=1,2,3$ as $(\mathbf{e}_{i,\parallel},\mathbf{e}_{i,\bot})$. All the geometrical parameters are gathered in Figure \ref{schema_nageur}).

Let us describe briefly  the forces applied to the robot.

\subsubsection{Elasticity}

The torsional springs which connect the swimmer segments exert a torque $\mathbf{T}^{\mathrm{el}}$ proportional to the shape angles $\alpha_1$ and $\alpha_2$. Thus, the torque $\mathbf{T}^{\mathrm{el}}_2$ exerted on $S_2$ is given by $\mathbf{T}^{\mathrm{el}}_2 = \kappa \alpha_1 \mathbf{e}_z$ and the torque $\mathbf{T}^{\mathrm{el}}_3$ exerted on $S_3$ is given by $\mathbf{T}^{\mathrm{el}}_3 = \kappa (\alpha_2-\alpha_0) \mathbf{e}_z$, where with $\alpha_0 \in (-\pi,\pi)$.

\begin{figure}
\begin{center}
\includegraphics[width=3.2in]{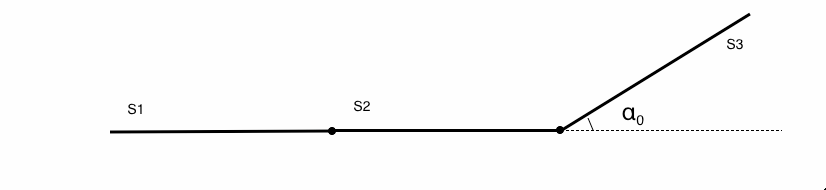}
\caption{The bent swimmer at its equilibrium.}
\label{equilibre}
\end{center}
\end{figure}

Here, the specificity of the swimmer we study is that the spring that relies $S_2$ and $S_3$ is at rest when $\alpha_2 = \alpha_0$. Hence, if $\alpha\not =0$, the springs tend to get the swimmer back to a bent shape, in which $S_1$, $S_2$ are aligned and the angle between $S_2$ and $S_3$ is equal to $\alpha_0$ (see Figure \ref{equilibre}). 

\subsubsection{Hydrodynamics}

Since the swimmer is assumed to be immersed in a fluid, hydrodynamic forces and torques derive from their interaction. According to the Resistive Force Theory (see \cite{GrayHancock55}), we assume that the drag force per unit length intensity is proportional to the velocity and to the hydrodynamics coefficients $\xi$ and $\eta$. Let $\mathbf{x}_s$ be a point on one of the segments $S_i$. Its velocity $\mathbf{u}_i (s)$ is given in the moving frame $(\mathbf{e}_{i,\parallel},\mathbf{e}_{i,\bot})$ by $\mathbf{u}_i (\mathbf{x}_s)=u_{i,\parallel} \mathbf{e}_{i, \parallel} + u_{i,\bot} \mathbf{e}_{i,\bot}$. The drag force exerted on this point is then given by 
$$
\mathbf{f}_i (\mathbf{x}_s) = -\xi_i u_{i,\parallel} \mathbf{e}_{i, \parallel} - \eta_i u_{i,\bot} \mathbf{e}_{i,\bot}.
$$
Let us integrate to obtain the total force $\mathbf{F}^h_i$ exerted on $S_i$ :
$$
\mathbf{F}^h_i = \int_{S_i} \mathbf{f}_i (\mathbf{x}_s) d\mathbf{x}_s.
$$
Moreover, given a point $\mathbf{x}_0$, the drag torque for $S_i$ with respect to $\mathbf{x}_0$ takes the form 
$$
\mathbf{T}^h_{i,\mathbf{x}_0} = \int_{S_i} (\mathbf{x}_s -\mathbf{x}_0 ) \times \mathbf{f}_i (\mathbf{x}_s) d\mathbf{x}_s.
$$
Hydrodynamic drag effects are resistant : they oppose to the swimmer's movement. Then, without a magnetic field, the swimmer tends to its equilibrium bent shape.

\subsubsection{Magnetism}

We assume that we apply a uniform time-varying external magnetic field $\mathbf{H}(t)$ in the fluid around the swimmer. Here, this magnetic field is assumed to be the control function and in the following, we express $\mathbf{H}(t)$ such that the robot can move along a prescribed trajectory. 
We choose to decompose $\mathbf{H}$ in the moving frame associated to $S_1$ : $\mathbf{H}(t)=H_{\parallel} \mathbf{e}_{1,\parallel} + H_{\bot} \mathbf{e}_{1,\bot}$. The magnetic field exerts a torque $\mathbf{T}^m_i$ on $S_i$ which is proportional to its magnetization coefficient $M_i$ : $\mathbf{T}^m_i = M_i \mathbf{e}_{i,\parallel} \times \mathbf{H}$.

\subsubsection{Dynamics equations}
The swimmer is considered sufficiently small to be at low Reynolds number regime, so that inertia may be neglected (see \cite{purcell1977life} for further considerations on low Reynolds number swimming). We apply Newton's second law to the system $\{ S_1 + S_2 + S_3 \}$ : the total force applied to the system is zero, and so is the total torque with respect to $\mathbf{x}$. Same holds for the subsystems $\{ S_2 + S_3 \}$ and $\{ S_3 \}$, with torques computed respectively with respect to $S_2$ end and $S_3$ end. It gives the following system of equations :
\begin{equation}
\left \{
\begin{array}{l l l c c}
\scriptstyle
\mathbf{F}_1^h + \mathbf{F}_2^h + \mathbf{F}_3^h & & & \scriptstyle = &\scriptstyle 0 \\
\scriptstyle
\mathbf{T}_{1,\mathbf{x}}^h + \mathbf{T}_{2,\mathbf{x}}^h+\mathbf{T}_{3,\mathbf{x}}^h & \scriptstyle + \mathbf{T}_1^m + \mathbf{T}_2^m + \mathbf{T}_3^m & & \scriptstyle = & \scriptstyle 0 \\
\scriptstyle
\mathbf{T}_{2,\mathbf{x}_2}^h + \mathbf{T}_{3,\mathbf{x}_2}^h & \scriptstyle + \mathbf{T}_2^m+\mathbf{T}_3^m &\scriptstyle + \quad \; \,\mathbf{T}_2^\mathrm{el} & \scriptstyle = & \scriptstyle 0 \\ 
\scriptstyle
\underbrace{\scriptstyle \mathbf{T}_{3,\mathbf{x}_3}^h \qquad \qquad \:}_{\mathrm{hydrodynamic \: terms}} &\scriptstyle + \underbrace{\scriptstyle \mathbf{T}_3^m \qquad \quad \:}_{\mathrm{magnetic \: terms}} &\scriptstyle + \underbrace{\scriptstyle \mathbf{T}_3^\mathrm{el}}_{\mathrm{elastic \: terms}} & \scriptstyle = & \scriptstyle 0
\end{array}
\right.
\end{equation}
This system gives five scalar equations by projecting the first line on $(Ox)$ and $(Oy)$ and the last three on $(Oz)$. After computing the different contributions, the system takes the form
\begin{equation}
M(\alpha_1,\alpha_2) R_{-\theta} \dot{Z} = Y,
\label{equation_de_base}
\end{equation}
with $Z=\begin{pmatrix} x & y &\theta & \alpha_1 & \alpha_2 \end{pmatrix} ^T$,
$$
R_{\theta} = \left (  \begin{array}{c | c} r_{-\theta} =\begin{pmatrix} \cos{\theta} & \sin{\theta} \\ -\sin{\theta} & \cos{\theta} \end{pmatrix} & 0 \\ \hline 0 & I_3 \end{array} \right ) 
$$
and
$$
Y= \left ( \begin{array}{c} 0 \\ 0 \\  \! \! \! \! \! \! \! \! \! \! \! \! \! \! \! \! \! \! \! \! \! \! \! \!\scriptstyle{H_{\parallel} (M_2 \sin{\alpha_1} + M_3 \sin{(\alpha_1 + \alpha_2)})} \\ \qquad \qquad \scriptstyle{- H_{\bot} (M_1 + M_2 \cos{\alpha_1} + M_3 \cos{(\alpha_1 + \alpha_2)})} \\ \! \! \! \! \! \! \! \! \! \! \! \! \scriptstyle{ -\kappa \alpha_1 + H_{\parallel} (M_2 \sin{\alpha_1} + M_3 \sin{(\alpha_1 + \alpha_2)})} \\ \qquad \qquad \scriptstyle{ - H_{\bot} ( M_2 \cos{\alpha_1} + M_3 \cos{(\alpha_1 + \alpha_2)})} \\ \scriptstyle{-\kappa (\alpha_2-\alpha_0) + H_{\parallel} M_3 \sin{(\alpha_1 + \alpha_2)} - H_{\bot}  M_3 \cos{(\alpha_1 + \alpha_2)}} \end{array} \right ).
$$
$M$ is a matrix that depends only on $\alpha_1$ and $\alpha_2$.


\begin{rmk}Up to a rotation matrix that can be eliminated by a changing of basis, the dynamics only depends on the shape state variables $\alpha_1$ and $\alpha_2$ : the problem is invariant by any translation or rotation.
\label{invariance}
\end{rmk}
\begin{rmk}
If the magnetic field is supposed to be zero, one can see that the equilibrium points are of the form $(x,y,\theta,0,\alpha_0)$ with $(x,y,\theta) \in \mathbf{R}^3$.
\end{rmk}

Straightforward computations show that the determinant of $M$ remains negative for all $(\alpha_1,\alpha_2)$, so $M$ is invertible and we can rewrite the system (\ref{equation_de_base}) as a nonlinear control system given by
\begin{equation}
R_{-\theta} \dot{Z} = \mathbf{F}_0 + H_{\parallel} (t) \mathbf{F}_1+ H_{\bot} (t) \mathbf{F}_2 ,
\label{systeme_controle_robot}
\end{equation}
where $\mathbf{F}_0$,$\mathbf{F}_1$ and $\mathbf{F}_2$ are combinations of the third, fourth and fifth columns of $M^{-1}$, denoted respectively in what follows by $\mathbf{X}_3$,$\mathbf{X}_4$ and $\mathbf{X}_5$ :
$$
\begin{array}{r l}
\mathbf{F}_0 =& -\kappa( \alpha_1 \mathbf{X}_4 + (\alpha_2-\alpha_0) \mathbf{X}_5); \\
\mathbf{F}_1=& (M_2 \sin{\alpha_1} + M_3 \sin{(\alpha_1+\alpha_2)})(\mathbf{X}_3 + \mathbf{X}_4) \\ & +M_3 \sin{(\alpha_1+\alpha_2)} \mathbf{X}_5 ;\\
\mathbf{F}_2 = & -M_1 \mathbf{X}_3 \\ & - (M_2 \cos{\alpha_1} + M_3 \cos{(\alpha_1 + \alpha_2)}) (\mathbf{X}_3 + \mathbf{X}_4)  \\ & - M_3 \cos{(\alpha_1+ \alpha_2)} \mathbf{X}_5.
\end{array}
$$

\subsection{Partial Controllability and Small-Time Local Controllability}
Let us remind that our aim is not to control either only the position of the swimmer or both its position and its orientation, without taking care of its shape. This type of problem is a \emph{partial controllability} problem, or \emph{$\Pi_p$-controllability} problem (see \cite{duprez2015controlabilite}). Let us define $\Pi_p$ as the projection operator given by $$
\begin{array}{c c c c}
\Pi_p : & \mathbf{R}^p \times \mathbf{R}^{n-p} & \rightarrow & \mathbf{R}^p \\
& (y_1,y_2) & \mapsto & y_1,
\end{array}
$$
where $n$ is the state space dimension for the considered control system and $1 \leq p \leq n$. Let us define the  \emph{$\Pi_p$-controllability} and state the classical Kalman condition for linear systems. 

\begin{mydef}
Let $(S)$ be the linear system of ordinary differential equations
$$
(S): \left \{ \begin{array}{l} \dot{y}=Ay+Bu \quad \text{in } [0,T] \\ y(0)=y_0, \end{array} \right . 
$$
with $y_0 \in \mathbf{R}^n$, $A \in M_n(\mathbf{R})$, $B \in M_{n,m} (\mathbf{R})$, and $u \in L^2([0,T],\mathbf{R}^m)$ the control. $(S)$ is $\Pi_p$-controllable at time $T$ if for all $y_0 \in \mathbf{R}^n$ and $y_T \in \mathbf{R}^p$, there exists $u$ such that the solution of (S) verifies 
$$
\Pi_p y(T;y_0,u)=y_T.
$$
\end{mydef}

\begin{mythm}\label{th1}
The system $(S)$ is $\Pi_p$-controllable at time $T$ if and only if 
$$
\text{Ker} (K^T \Pi_p^T) = \{ 0 \},
$$
where $K$ is the \emph{Kalman matrix} given by $K=\begin{pmatrix} B & AB & A^2 B & \dots & A^{n-1} B \end{pmatrix}$. 
\label{kalman}
\end{mythm}

In other terms, $(S)$ is $\Pi_p$-controllable if and only if the submatrix of K consisting of the $p$ first rows of $K$ is of maximal rank $p$. 

For nonlinear systems, global controllability results such as above are often hard to obtain. Hence, we aim for local results, such as \emph{small-time local partial controllability} (abbreviated as STLPC). Let $(NL)$ be the linear system of ordinary differential equations
$$
(NL):(\dot{y}_1,\dot{y}_2) = f(y_1,y_2,u),
$$
where  $y_1 \in L^2([0,T],\mathbf{R}^p)$, $u \in L^2([0,T],\mathbf{R}^{n-p})$, $u \in L^2([0,T],\mathbf{R}^m)$ and $f: \mathbf{R}^p \times \mathbf{R}^{n-p} \times \mathbf{R}^m$ verifies the condition of the global Cauchy-Lipschitz Theorem.
\begin{mydef}
Let $(y_1^e,y_2^e,u^e) \in \mathbf{R}^p \times \mathbf{R}^{n-p} \times \mathbf{R}^m$ be an equilibrium of a control system (NL). The control system  (NL) is \emph{small-time locally partially controllable at $(y_1^e,y_2^e,u^e)$}  with respect to  $y_1$ if for every $\epsilon >0$, there exists a real number $\eta >0$ such that, for every $(y_{1}^0,y_2^0,y_{1}^f) \in B_{\eta}(y_1^e) \times B_{\eta}(y_2^e) \times B_{\eta}(y_1^e)$, there exists $u\in L^2([0,\epsilon] \rightarrow \mathbf{R}^m)$ such that $(y_1,y_2)$ verifies system (NL),
$$
\begin{array}{l c}
\mathrm{(i)} & \forall t \in [0,\epsilon], | u(t)-u^e | \leq \epsilon; \\
\mathrm{(ii)} & y_1(\epsilon)=y_1^f.
\end{array}
$$
\end{mydef}
An immediate application of the inverse mapping theorem enables us to obtain the following usual sufficient condition for the STLPC around an equilibrium for a nonlinear system.

\begin{mythm}
The nonlinear control system (NL)  is STLPC at an equilibrium if its linearized control system around this equilibrium is $\Pi_p$-controllable for some time $T>0$.
\label{stplc}
\end{mythm}

This last theorem justifies that we study the linearized system of \eqref{systeme_controle_robot} to get partial controllability around the equilibrium.

\subsection{Local partial controllability result}

In the following, we prove that the position of the 3-link magnetic swimmer can be partially controlled by the external magnetic fields when the bent swimmer is close to its equilibrium. The main result states as follows.

\begin{mythm}
If $\alpha_0 \neq 0$, then system \eqref{systeme_controle_robot} is STPLC with respect to $(x,y)$ around any equilibrium point.
\label{result}
\end{mythm}

\emph{Proof}. We only need to prove the result for the particular equilibrium point $\mathbf{O}=(0,0,0,0,\alpha_0)$. Indeed, according to Remark \ref{invariance}, solutions of \eqref{systeme_controle_robot} are invariant under the transformations 
$$ \begin{array}{l}
\left ( \begin{pmatrix} x \\ y \end{pmatrix},\theta,\alpha_1,\alpha_2,H_{\parallel},H_{\bot} \right ) \\ \quad \quad \mapsto \left ( R_{\bar{\theta}} \begin{pmatrix} x+\bar{x} \\ y + \bar{y} \end{pmatrix}, \theta + \bar{\theta},\alpha_1,\alpha_2,H_{\parallel},H_{\bot} \right ), \end{array}$$
 so if the result holds for $\mathbf{O}$, it may be carried to an arbitrary equilibrium point. 
Following Theorem \ref{stplc}, we look at the linearized system around the equilibrium point $\mathbf{O}$, which is given by
\begin{equation}
\dot{Z}=A Z+B \mathbf{H},
\label{linearsystem}
\end{equation}
where $A$ is the Jacobian of $Z \mapsto \mathbf{F}_0 (Z)$ at $\mathbf{O}$ and the two columns of $B$ are $\mathbf{F}_1 (\mathbf{O})$ and $\mathbf{F}_2 (\mathbf{O})$. Since we are interested in controlling on the position $(x,y)$,
we are looking for the $\Pi_2$ partial controllability  of the swimmer, hence, according to Theorem \ref{th1}, we have to check that the first two rows of the Kalman matrix $$K=\begin{pmatrix} B & AB & A^2 B & A^3 B & A^4 B \end{pmatrix}$$ give a matrix of rank 2. If we look at the $2\times 2$ submatrix given by \footnote{It would be more natural to chose the two columns of $B$. Indeed, for \emph{almost all} values of the parameters, the two first lines of $B$ compose an invertible $2\times2$ matrix (its determinant is called $D(0,\alpha_0)$ in Section \ref{explicit}). For a precise set of parameters, however, this matrix has rank 1 only and one has to choose the first columns of $B$ and $AB$ instead; since this choice works in general, it is the one we make in the proof.} the first two entries of the first columns of $B$ and $AB$, a tedious but straightforward computation enables us to get the determinant of this matrix :
\begin{equation}
\textstyle{\frac{108 M_3^2 \kappa ( -9 \eta \xi (19 \eta + 54 \xi) \cos \alpha_0 - 2 \Xi (\eta + 2 \xi) ) \sin^3 (\alpha_0)}{L^7 \eta^2 (\eta^2 + 34 \eta \xi + 28 \xi^2 - (\eta^2-11 \eta \xi + 28 \xi^2) \cos (2\alpha_0))^2},}
\label{determinant_sub}
\end{equation}
where 
$$
\textstyle{\Xi = (\eta^2 + 19 \eta \xi + 7 \xi^2 - (\eta^2 - 8 \eta \xi + 7 \xi ^2) \cos (2 \alpha_0) ).}
$$
Along with the hypothesis $\alpha_0 \neq 0$, the straightforward inequalities 
$$
\eta^2 + 19 \eta \xi + 7 \xi^2 > (\eta^2 - 8 \eta \xi + 7 \xi ^2) \cos (2 \alpha_0)
$$
and
$$
\eta^2 + 34 \eta \xi + 28 \xi^2 > (\eta^2-11 \eta \xi + 28 \xi^2) \cos (2\alpha_0)
$$
show that the numerator and denominator of (\ref{determinant_sub}) are respectively the opposite sign of $\alpha_0$ and positive. Therefore, the determinant is nonzero and its associated submatrix has rank 2. According to Theorem \ref{kalman}, the system (\ref{linearsystem}) is partially controllable. 
We conclude by applying Theorem \ref{stplc}. 

\begin{rmk} 
\emph{None of the above applies for a non-bent 3-link swimmer, i.e. if $\alpha_0 = 0$.}
Indeed, the proof does not work in this case because the numerator of the determinant \eqref{determinant_sub} is zero.
Furthermore, a straightforward computation yields
$$
A=\begin{pmatrix} 0 & 0 & 0 & 0 & 0 \\ 0 & 0 & 0 & * & * \\ 0 & 0 & 0 & * & * \\ 0 & 0 & 0 & * & * \\ 0 & 0 & 0 & * & * \end{pmatrix} \quad \text{and} \quad B=\begin{pmatrix} 0 & 0 \\ 0 & * \\ 0 & * \\ 0 & * \\ 0 & * \end{pmatrix},
$$
where stars stand for possibly nonzero entries,
hence the first row of
the 
Kalman matrix is zero, the linearized control system is not controllable, and Theorem \ref{stplc} does not apply.

We do not know whether \eqref{systeme_controle_robot} is STPLC around its equilibrium if $\alpha_0 = 0$,
but non controllability of the linearized system in this case indicates that the non-bent 3-link swimmer is harder to control. It is why we focus here on the case of bend swimmer.
\end{rmk}

\section{EXPLICIT PARTIAL CONTROL \label{explicit}}

%

In this Section, we describe a method to make the swimmer's position $(x,y)$ follow an arbitrary trajectory, while its orientation and its shape are not prescribed, and we explain the theoretical difficulties that arise.

\begin{figure}
\begin{center}
 \includegraphics[width=3.2in]{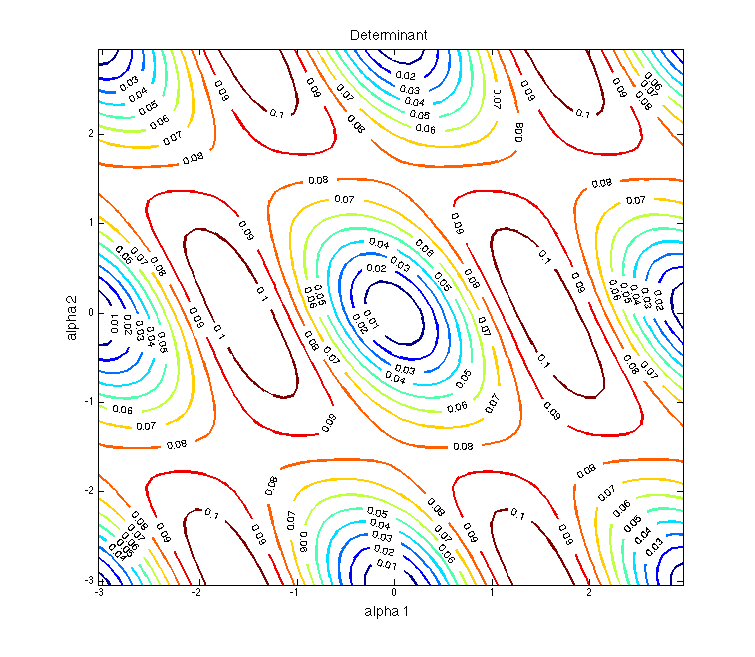}
\caption{Aspect of $D(\alpha_1,\alpha_2)$ with numerical values from Table \ref{table_values}. It vanishes only at $(0,0)$.}
\label{deter}
\end{center}
\end{figure}

Let us focus on the two first lines of the system \eqref{systeme_controle_robot}:
\begin{equation}
r_{-\theta} \begin{pmatrix} \dot{x} \\ \dot{y} \end{pmatrix} = \begin{pmatrix} F_{0x} + H_{\parallel} F_{1x} + H_{\bot} F_{2x} \\ F_{0y} + H_{\parallel} F_{1y} + H_{\bot} F_{2y} \end{pmatrix} \,.
\label{equations_xy}
\end{equation}
We denote by $D(\alpha_1,\alpha_2)$ the determinant of the $2\times2$ matrix $\begin{pmatrix} F_{1x} & F_{2x} \\ F_{1y} & F_{2y} \end{pmatrix}$; it depends only on the state variables $\alpha_1,\alpha_2$. It is clear that $D$ vanishes at straight positions ($D(0,0)=0$). Restricting $D$ to the open square $K=(-\pi,\pi) \times (-\pi,\pi)$
(values of $(\alpha_1,\alpha_2)$ outside $K$ are not physical: the segments would then overlap),
and using the values of the parameters used for our numerical simulations (see Table \ref{table_values}),
the plot in Figure \ref{deter} shows that it vanishes only at $(0,0)$. 

Let $T>0$. Let $f$ and $g$ be two functions of class $C^1$ on $[0,T]$. We require that the swimmer follows exactly the trajectory parameterized by $f$ and $g$, i.e.
\begin{equation}
\forall t \in [0,T], \begin{pmatrix} x(t) \\ y(t) \end{pmatrix} = \begin{pmatrix} f(t) \\ g(t) \end{pmatrix}\,.
\label{to_diff}
\end{equation}
The problem is to find the control functions $H_{\parallel},H_{\bot}$ that achieve this goal.  Differentiating \eqref{to_diff} and using (\ref{equations_xy}), we get 
$$
\begin{pmatrix} F_{0x} + H_{\parallel} F_{1x} + H_{\bot} F_{2x} \\ F_{0y} + H_{\parallel} F_{1y} + H_{\bot} F_{2y} \end{pmatrix} = r_{\theta} \begin{pmatrix} f'(t) \\ g'(t) \end{pmatrix}.
$$
Hence, at each time $t$, $H_{\bot}$ and $H_{\parallel}$ must solve the $2\times2$ linear system of equations 
\begin{equation}
\begin{pmatrix} F_{1x} & F_{2x} \\ F_{1y} & F_{2y} \end{pmatrix} \begin{pmatrix} H_{\parallel} \\ H_{\bot} \end{pmatrix} = \begin{pmatrix} -F_{0x} \\ - F_{0y} \end{pmatrix} + r_{\theta} \begin{pmatrix} f'(t) \\ g'(t) \end{pmatrix}.
\label{controls_trajectory}
\end{equation}
It has a unique solution if the determinant $D(\alpha_1(t),\alpha_2(t))$ does not vanish.

The functions $H_{\parallel}$ and $H_{\bot}$ solving the system \eqref{controls_trajectory} depend on $f'$, $g'$ and the state variables $\theta$, $\alpha_1$, $\alpha_2$, $x$ and $y$.
Following \cite[Chapter 7]{Isid95}, this can be re-formulated in terms of "relative degree" and non-interactive control. Considering the control system \eqref{systeme_controle_robot} with two outputs $f=x$ and $g=y$, and inputs $H_{\parallel},H_{\bot}$, it has vector relative degree $\{1,1\}$ when $D\neq0$ and the above mentioned expression of $H_{\parallel},H_{\bot}$ as functions of $f'$, $g'$, $\theta$, $\alpha_1$, $\alpha_2$, $x$, $y$ is a feedback transformation that solves the ``Noninteracting Control Problem'' \cite[Section 7.3]{Isid95}, i.e. it produces a control system with new controls $f'$ and $g'$ where the control $f'$ acts on the output $x$ only and the control $g'$ acts on the output $y$ only, as can be seen by computing $\dot x=f'$ and $\dot x=g'$. Clearly, for any functions $f(t), g(t)$, if such a feedback transformation exists and if $f'(t),g'(t)$ are taken to be the time derivatives of $f(t),g(t)$ and if $x(0)=f(0)$ and $y(0)=g(0)$, then \eqref{equations_xy} holds.
Note that, because of the intermediary state feedback, the expression of $H_{\parallel},H_{\bot}$ \emph{as a function of time} is obtained implicitly after solving the closed-loop ODE.

This method \emph{fails} if the swimmer has to go through the straight configuration, where one cannot invert~\eqref{controls_trajectory}. 
Since we cannot guarantee in general that the straight confuguration will not be encountered, all we can do is try numerically, see next section where we see that both may happen:
either the method works until the end of the trajectory or the alignment occurs and the controls blow up.
This does not contradict Theorem~\ref{result} that is local and only concerned with short trajectories; Theorem~\ref{result} is however stronger since it applies even if $D$ vanishes at the equilibrium (see the footnote on the thuird page).

%
%
%

\section{NUMERICAL SIMULATIONS \label{numericalsimu}}
In this part, we show numerical result obtained by using the approach given in the previous Section, in the case where $D$ vanishes only at $\alpha_1=\alpha_2=0$. 
Table \ref{table_values} gives the actual values for the parameters used in this Section (the values have been chosen according to \cite{alouges2013self}). We have used Matlab and more particularly the function \textit{ode23t} to integrate the system \eqref{systeme_controle_robot}.

\begin{figure}
\begin{center}
\includegraphics[width=2.5in]{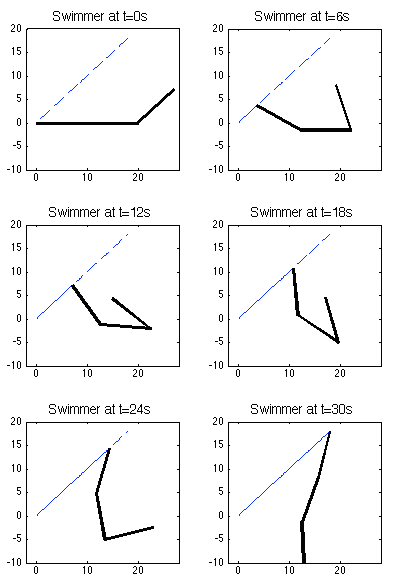}
\caption{Snapshots of the swimmer (in black) following a straight line (in blue). The blue plain line indicates the path already done, and the blue dotted line indicates the remaining path. The scale is in micrometers.}
\label{straight_line}
\end{center}
\end{figure}

\begin{figure}
\begin{center}
\includegraphics[width=1.6in]{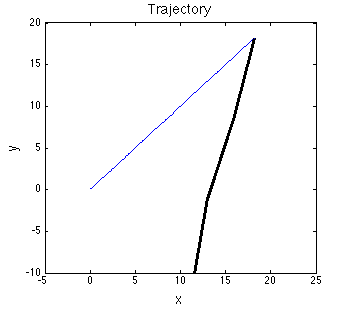}
\caption{An example of a bad case in which the swimmer is aligning.}
\label{bad_traj}
\end{center}
\end{figure}

\begin{figure}
\begin{center}
\includegraphics[width=2.4in]{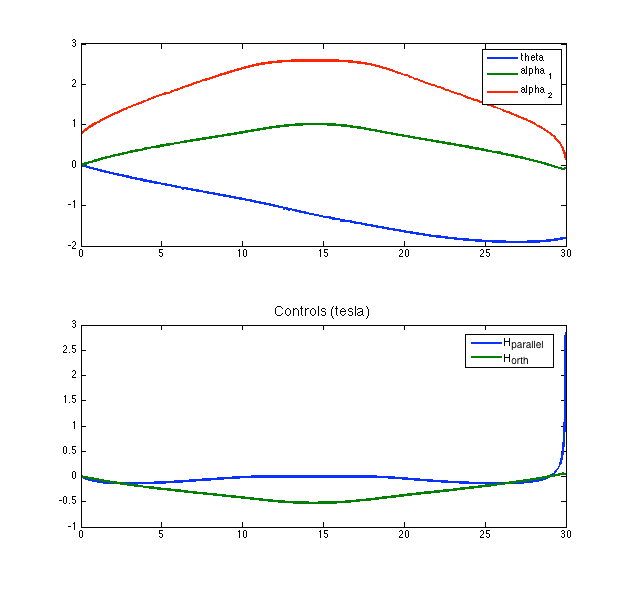}
\caption{State variables and controls along the straight line. At the end of the time, we can see that $\alpha_1$ and $\alpha_2$ are going to zero and that $H_\parallel$ goes to infinity.}
\label{bad_controls}
\end{center}
\end{figure}

The routine takes functions $f$ and $g$ as an input, and returns the trajectory of the robot, including its orientation and shape, and the required controls $H_{\parallel}$ and $H_{\bot}$. 
Figure \ref{straight_line} shows how it can follow a straight line. 
When the swimmer is closed to the aligned position, according to the system \eqref{controls_trajectory}, the magnetic field goes to infinity, then the simulation stops and we cannot follow the entire prescribed trajectory. It is the case in Figure \ref{bad_traj} and Figure \ref{bad_controls}. 

\begin{figure}
\begin{center}
\includegraphics[width=2.5in]{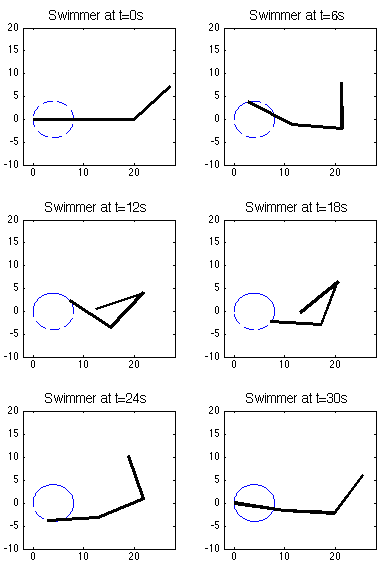}
\caption{Snapshots of the swimmer (in black) following a circle (in blue) and going back to its initial position. The blue plain line indicates the path already done, and the blue dotted line indicates the remaining path.}
\label{cercle}
\end{center}
\end{figure}

\begin{figure}
\begin{center}
\includegraphics[width=3in]{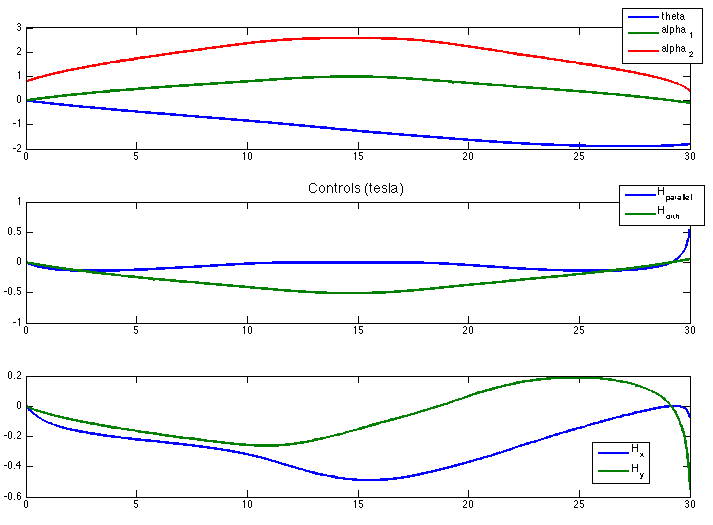}
\caption{State variables and controls along the circle trajectory.}
\label{controles_cercle}
\end{center}
\end{figure}

\begin{figure}
\begin{center}
\includegraphics[width=3.2in]{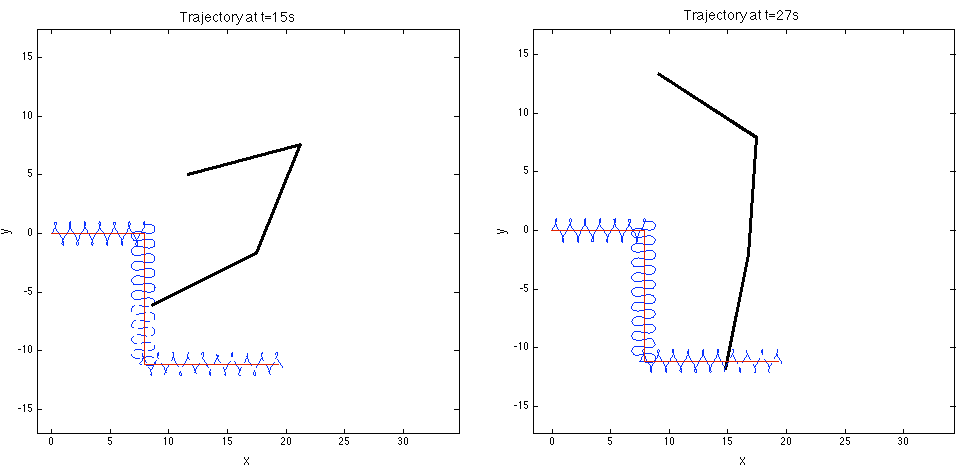}
\caption{Snapshots of the swimmer (in black) following a complex trajectory (in blue). The blue plain line indicates the path already done, and the blue dotted line indicates the remaining path. The red line indicates the global direction that we want to follow.}
\label{trajectory1}
\end{center}
\end{figure}

\begin{figure}
\begin{center}
\includegraphics[width=2.8in]{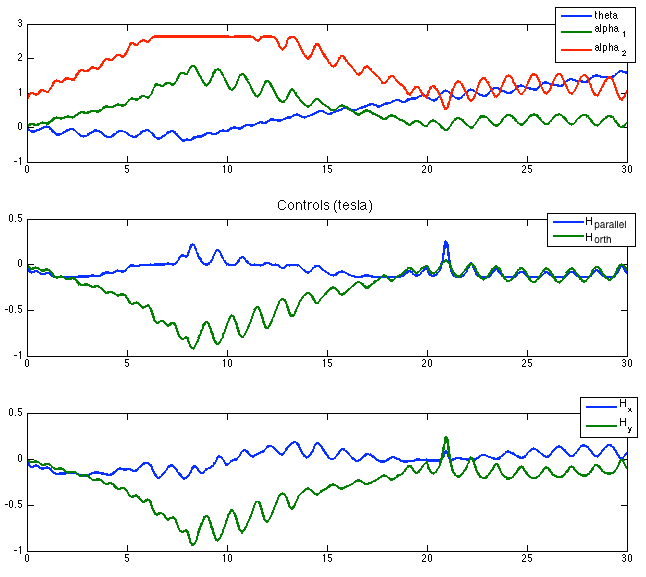}
\caption{Angles and controls with respect the time (in sec.) during the trajectory described in Fig. \ref{trajectory1}. The middle plot shows the components $H_{\parallel}$ and $H_{\bot}$ of the control, whereas the bottom plot shows the $H_{x}$ and $H_{y}$ components, in the reference basis.}
\label{controlsangles1}
\end{center}
\end{figure}

Figure \ref{cercle} shows a small circle trajectory that the swimmer is able to follow. Figure \ref{controles_cercle}  presents its orientation, its shape and the external magnetic fields along this experiment. 
Finally, Figure \ref{trajectory1} shows a more complex trajectory that leads to follow a path while remaining close to a global direction.
Figure \ref{controlsangles1} shows the associated controls and angular variables during it.

\begin{table}
\begin{center}
\begin{tabular}{| l | l |}
\hline
Parameter & Value \\
\hline
$\ell$ & 10 $\mu$m \\
$\eta$ & 12.4$\times 10^{-3}$ N.s.m$^{-2}$ \\
$\xi$ & 6.2$\times 10^{-3}$ N.s.m$^{-2}$ \\
$M_1$ & 1.6 A.$\mu \text{m}^2$ \\
$M_2$ & 2.4 A.$\mu \text{m}^2$ \\
$M_3$ & 3.2 A.$\mu \text{m}^2$ \\
$\kappa$ & 8.3$\times 10^{-7}$ N.$\mu$m \\
\hline
\end{tabular}
\caption{Numerical values used for the parameters}
\label{table_values}
\end{center}
\end{table}

%
%

\section{CONCLUSION AND PERSECTIVES}
\label{perspectives}

Here, we proved local partial controllability of the ``bent'' swimmer.
The same result for the \emph{non bent} 3-link swimmer is still an open question, the present work indeed stemmed out from trying  to prove STPLC via the return method of Coron~\cite[Chapter 6]{coron2007control}, rather than Theorem \ref{stplc} (linear test).

 Trying to go beyond local for the bent swimmer, we described a method to drive the position of the swimmer along a given trajectory.
It fails if the swimmer passes through the straight shape, and we cannot ensure that this will not occur.
In our numerical experiments, we observe that the method applies to some trajectories but that in other cases one has to go through the straight shape and the controls blow up, evidencing that this straight shape still represents a serious barrier to maneuverability, even though bending the swimmer provides linear controllability at the equilibrium. 

Another perspective, under our investigation, is to conduct further numerical study, and add energetic aspects, in order to find a magnetic field which allows the swimmer to move close to a prescribed path by minimizing the kinetic energy of fluid-swimmer system (see \cite{Lighthill76}). 



\end{document}